\documentclass[11pt]{article}
\usepackage[psamsfonts]{amsfonts}
\usepackage{amsmath, amssymb, amsthm}

\newcommand{\bt}   {\begin{theorem}}
\newcommand{\et}   {\end  {theorem}}
\newcommand{\bl}   {\begin{lemma}}
\newcommand{\el}   {\end  {lemma}}
\newcommand{\bp}   {\begin{prop}}
\newcommand{\ep}   {\end  {prop}}
\newcommand{\bc}   {\begin{cor}}
\newcommand{\ec}   {\end  {cor}}
\newcommand{\bd}   {\begin{defn}}
\newcommand{\ed}   {\end  {defn}}

\newcommand{\ba}   {\begin{array}}
\newcommand{\ea}   {\end  {array}}
\newcommand{\be}   {\begin{enumerate}}
\newcommand{\ee}   {\end  {enumerate}}
\newcommand{\bi}   {\begin{itemize}}
\newcommand{\ei}   {\end  {itemize}}

\def\eq#1\en{\begin{equation}#1\end{equation}}  
\def\eqsplit#1\ensplit{
	\begin{equation}\begin{split}#1\end{split}\end{equation}
	}
\def\eqalign#1\enalign{
	\begin{align}#1\end{align}
	}
\def\eqmul#1\enmul{
	\begin{multline}#1\end{multline}
	}
\newcommand{\eqarrstar} {\begin{eqnarray*}} 
\newcommand{\enarrstar} {\end{eqnarray*}} 
\newcommand{\eqarray}   {\begin{eqnarray}} 
\newcommand{\enarray}   {\end{eqnarray}}

\newcommand{\lbeq}[1]  {\label{e:#1}}
\newcommand{\refeq}[1] {\eqref{e:#1}}    

\makeatletter
\newcommand{\labelcounter}[2]{{%
	\stepcounter{#1}
	\protected@write\@auxout{}%
	{\string\newlabel{#2}{{\csname the#1\endcsname}{\thepage}}}%
	{\ref{#2}}
	}}
\makeatother
\newcommand{\sss}   { \scriptscriptstyle } 

\newcommand{\spose}[1] {{\hbox to 0pt{#1\hss}} }
\newcommand{\ltapprox} {\mathrel{\spose{\lower 3pt\hbox{$\mathchar"218$}}
 \raise 2.0pt\hbox{$\mathchar"13C$}}}
\newcommand{\gtapprox} {\mathrel{\spose{\lower 3pt\hbox{$\mathchar"218$}}
 \raise 2.0pt\hbox{$\mathchar"13E$}}}

\numberwithin{equation}{section}

\newcommand{\een}{\end{enumerate}}
\newcommand{\noi}{\noindent}
\newcommand{\ra}{\rightarrow}
\newcommand{\sig}{\sigma}
\newcommand{\eps}{\varepsilon}
\newcommand{\del}{\delta}

\newcommand{\TREE}{{\rm TREE}}
\newcommand{\NICE}{{\rm NICE}}
\newcommand{\BIN}{{\rm BIN}}
\newcommand{\MIN}{{\rm MIN}}

\newcommand{\ESC}{{\rm ESC}}
\newcommand{\FAIL}{{\rm FAIL}}

\newcommand{\normal}{{\rm normal}}
\newcommand{\Po}{{\rm Po}}
\newcommand{\E}{E}
\renewcommand{\P}{\Pr}

\newcommand{\lam}{\lambda}
\newcommand{\vep}{\varepsilon}
\newcommand{\ah}{\alpha}
\newcommand{\nn}{\nonumber}

\newcommand{\ben}{\begin{enumerate}}
\newcommand{\beq}{\begin{equation}}
\newcommand{\eeq}{\end{equation}}
\newcommand{\bec}{\begin{center}}
\newcommand{\ece}{\end{center}}

\newtheorem{definition}{Definition}

\newtheorem{theorem}{Theorem}[section]

\newtheorem{prop}[theorem]{Proposition}
\newtheorem{lemma}[theorem]{Lemma}

\newtheorem{cor}[theorem]{Corollary}

\newcommand{\N}{\mathbb{N}}

\newcommand{\sR}{{\sss R}}
\newcommand{\sL}{{\sss L}}
\newcommand{\sY}{{\sss Y}}
\newcommand{\sM}{{\sss M}}

\title{Counting Connected Graphs Asymptotically}

\author{
    Remco van der Hofstad\thanks{Department of Mathematics and
        Computer Science, Eindhoven University of Technology,
        5600 MB Eindhoven, The Netherlands.
        {\tt rhofstad@win.tue.nl}}\\
    Joel Spencer\thanks{Department of Computer Science,
Courant Institute of Mathematical Sciences,
New York University, 251 Mercer St., New York, NY 10012, U.S.A.
{\tt spencer@cims.nyu.edu}}
}

\begin{document}
\maketitle


\begin{abstract}
We find the asymptotic number of connected graphs with $k$ vertices and
$k-1+l$ edges when $k,l$ approach infinity, reproving a result of Bender,
Canfield and McKay.  We use the {\em probabilistic method}, analyzing
breadth-first search on the random graph $G(k,p)$ for an appropriate edge
probability $p$.  Central is analysis of a random walk with fixed beginning
and end which is tilted to the left.
\end{abstract}

\section{The Main Results}\label{section1}
In this paper, we investigate the number of graphs with a given complexity.
Here, the {\em complexity} of a graph is its number of
edges minus its number of vertices plus one.  For $k,l\geq 0$,
we write $C(k,l)$ for the number of labeled connected graphs with $k$ vertices
and complexity $l$.\\

The study of $C(k,l)$ has a long history.
Cayley's Theorem gives the exact formula for the
number of trees, $C(k,0)=k^{k-2}$.  The asymptotic formula for
the number of unicyclic graphs, $C(k,1)$, has been given by
R\'enyi \cite{renyi} and others. Wright \cite{wright} gave
the asymptotics of $C(k,l)$ for $l$ arbitrary but fixed and $k\ra \infty$,
and also studied the asymptotic behavior of $C(k,l)$ when $l=o(k^{1/3})$
in \cite{wrightb}.

The asymptotics of $C(k,l)$ for all
$k,l\ra \infty$ were found by Bender, Canfield, and McKay \cite{bcm}.
The proof in \cite{bcm} is based on an explicit recursive formula for
$C(k,l)$. In this paper, we give an alternate, and
substantially different, derivation of the Bender, Canfield, McKay
results.  Our argument is Erd\H{o}s Magic, using the study of the random
graph $G(k,p)$ to find the asymptotics of the strictly enumerative $C(k,l)$.
The critical idea, given in Theorem \ref{exactthm} below, involves an
analysis of a breadth first search algorithm on $G(k,p)$.
Similar methods, with somewhat weaker results in our
cases, were employed recently by Coja-Oghlan, Moore and Sanwalani \cite{CMS}.
We can also use the results and methodology in this paper
to find local statistics on the
joint distribution of the size and complexity of the dominant
component of $G(n,p)$ in the supercritical regime, which we defer
to a future publication \cite{HS05}. Further,
while computational issues are not addressed in our current work,
these methods may be used to efficiently generate a random connected
graph of given size and complexity \cite{AS05}.
The idea of using the random graph to study $C(k,l)$ has appeared
previously. In \cite{luczak}, a reformulation in terms of random
graphs was used to prove upper and lower bounds on $C(k,l)$, extending
the upper bound in \cite{Boll84}. The idea in \cite{luczak} is that the
expected value of the number of connected components with $k$ vertices
and complexity $l$ can be explicitly written in terms of $C(k,l)$. Bounds on
the random number of such components then imply bounds on $C(k,l)$.
In \cite{Spen97}, a more sophisticated analysis was performed, where
the connected component of a given node in the random graph was
explored using breadth-first search. This analysis allows to rewrite
the asymptotics of $C(k,l)$ for $l$ fixed in terms of
Brownian excursions. Interestingly, this identifies the Wright constants
for the asymptotics $C(k,l)\sim c_l k^{k-2} k^{3l/2}$ in terms of moments of
the mean distance from the origin of a Brownian excursion. These moments
were also investigated in \cite{Louc84}, but the connection to the
Wright constants had not been made before.

In this case, we will use the breadth-first search representation
of connected components in $G(k,p)$ for $C(k,l)$, choose $p$
appropriately, and analyse the resulting problem using probabilistic
means. The critical identity is Theorem \ref{exactthm}, which rewrites
the $C(k,l)$ in terms of a $k$-step conditioned random walk with steps that
are Poisson random variables minus one, and where the parameter of the Poisson
steps varies with time. The main work in this paper then lies in the study
of this random walk.

We note that when $l\geq k\ln k$ (and even somewhat less) the asymptotics
of $C(k,l)$ are trivial as asymptotically almost all graphs on $k$ vertices
and $k-1+l$ edges will be connected.  Thus, while our methods extend further,
we shall restrict ourselves to finding the asymptotics of $C(k,l)$ with
$l\leq k\ln k$ and $l\ra \infty$.  It will be convenient to subdivide the
possible $l$ into three regimes:
\ben
\item Very Large: $l\gg k$ and $l\leq k\ln k$;
\item Large: $l=\Theta(k)$;
\item Small: $l\ll k$ and $l\ra \infty$.
\een

The main ideas of this paper are given in Section \ref{section1},  where
we state the main results Theorems \ref{probtreethm} and \ref{main2}.
The consequences of our main results for $C(k,l)$ are formulated in Section
\ref{sec-asyC}.  The
proofs of Theorems \ref{probtreethm} and \ref{main2} are given in
Section \ref{section2}.

\subsection{Tilted Balls into Bins}
\label{statements}
Let $k\geq 2$ be an integer. In this section, we define a
process of placing $k-1$ balls into $k$ bins with a tilted
distribution, which makes it more likely that balls are
placed in the left bins. In Section \ref{theconnection} below,
we will find an identity between this bin process and
$C(k,l)$.

Let $p\in (0,1]$. We have $k-1$ balls $1\leq j\leq k-1$
and $k$ bins $1\leq i \leq k$. We place the $j^{\rm th}$ ball
into bin $T_j$, where $T_j$ is a random variable with distribution
given by
    \beq\label{dist}
    \Pr[T_j=i]= \frac{p(1-p)^{i-1}}{1-(1-p)^k}.
    \eeq
This is a truncated geometric distribution.  Note that the larger $p$ is the more
the $T_j$ are tilted to the left.  We shall call $p$ the {\em tilt} of
the distribution. The $T_j$ are independent and identically
distributed.  Let $Z_i$, $1\leq i\leq k$, denote
the number of balls in the $i^{\rm th}$ bin.  Set $Y_0=1$ and $Y_i=Y_{i-1}+Z_i-1$
for $1\leq i\leq k$.  Note that $Y_k=0$ as there are precisely $k-1$ balls.
Let $\TREE$ be the event that
    \beq
    \label{treeY}
    Y_t > 0 \quad \mbox{ for }\quad 1\leq t\leq k-1,
    \eeq
or, alternatively,
    \beq
    \label{tree}
    Z_1+\ldots+Z_t\geq t \quad \mbox{ for }\quad 1\leq t\leq k-1.
    \eeq
We note that we use the term $\TREE$ because there is a natural bijection between
placements of $k-1$ balls into $k$ bins satisfying $\TREE$ and trees
on $k$ vertices.  Alternatively, one may consider $Y_0,Y_1,\ldots,Y_k$,
with $Y_0=1, Y_k=0$,  as a walk with fixed endpoints, or a bridge.
The condition $\TREE$ can
then be interpreted as saying that the bridge is an {\it excursion}.  In certain
limiting situations the bridge approaches a biased Brownian bridge, where the bias
depends on the parameter $p$.

\begin{definition}
\label{def-M}
Set
    \beq\label{defM}
    M= {k\choose 2} - \sum_{j=1}^{k-1} T_j
    \eeq
in the probability space in which the $T_j$ are independent with
distribution given by (\ref{dist}).   Set $M^*$ equal the same
random variable but in the above probability space conditioned
on the event $\TREE$.
\end{definition}

We can give an alternative definition of $M$ as follows:
    \beq\label{twosameeq}
    {k\choose 2} - \sum_{j=1}^{k-1}T_j = \sum_{i=1}^{k-1} (Y_i-1),
    \eeq
which can be seen by noting that both sides of
(\ref{twosameeq}) increase by one when one ball
is moved one position to the left and decrease by one
when one ball is moved one position to the right.
Since one can get from any placement
to any other placement via a series of these moves, the two sides
of (\ref{twosameeq}) must differ by a constant.
However, when $T_j=j$ for $1\leq j\leq k-1$, we have
$Y_i=1$ for $1\leq i\leq k-1$ and $Y_k=0$ and so the sides
are equal for this placement of balls.

\subsection{The Critical Identity}\label{theconnection}

The main idea of our approach is given in Theorem \ref{exactthm} below.
Note that this result is exact, there are no asymptotics.

\begin{theorem}\label{exactthm}
For all $k,l\in \N$, $p\in (0,1]$,
    \beq\label{x} A_1A_2A_3 =
    C(k,l)p^{k+l-1}(1-p)^{{k\choose 2}-(k+l-1)},
    \eeq
where
    \begin{eqnarray}
    \label{A1}
    A_1&=&(1-(1-p)^k)^{k-1},\\
    \label{A2}
    A_2&=&\Pr[\TREE],\\
    \label{A3}
    A_3&=&\Pr[\BIN[M^*,p]=l].
    \end{eqnarray}
\end{theorem}

\proof
The right hand side of (\ref{x}) is the probability that $G(k,p)$
is connected and has complexity $l$. We show that the left hand side
of (\ref{x}) also gives this probability.
Designate a root
vertex $v$ and label the other vertices $1\leq j\leq k-1$.  We
analyze breadth-first search on $G$, starting with root $v$.
(More precisely, the queue is initially $\{v\}$.
In Stage 1 we pop $v$ off the queue
and add to the queue the neighbors of $v$. Each successive stage we pop
a vertex off the queue and add to the queue its neighbors that haven't
already been in the queue.  The process stops when the queue is empty.)

Each non-root $j$ flips a coin $k$ times, heads with probability $p$.
The $i^{\rm th}$ flip being heads means in the breadth-first stage
that if the $i^{\rm th}$ stage is reached and $j$ has not yet entered the queue,
then $j$ is adjacent to the ``popped" vertex.  To get all vertices
it is necessary that each $j$ has at least one head.  This happens
with probability $A_1$.  Conditioning on that, we let $T_j$ be that first
$i$ when $j$ had a head.  So $T_j$ has the truncated geometric distribution of
(\ref{dist}).  While the process continues $Y_t$ is the size of the
queue.  The condition that the process doesn't terminate before
stage $k$ is precisely that no $Y_t=0$ for $1\leq t\leq k-1$
which is $\TREE$, so this gives $A_2$.  Now the only $\{w_1,w_2\}$
whose adjacency has not been determined are those for which (letting $w_1$
be the first one popped) $w_2$ was in the queue when $w_1$ was popped.
There are precisely $\sum_{t=0}^{k-1}(Y_t-1)$ of such pairs, i.e., we add
the size of the queue minus the popped vertex over each stage, except
for the last stage.
Since we are conditioning on $\TREE$, the random variable
$\sum_{t=0}^{k-1}(Y_t-1)$ has distribution $M^*$.
We now look at those pairs, each is adjacent with independent probability $p$
and to have complexity $l$, we need to have exactly $l$ such pairs
adjacent, so that the probability of this event equals $A_3$.
\qed
\vskip 0.5cm

\noindent
Our approach to finding the asymptotics of $C(k,l)$ will be to find
the asymptotics of $A_2,A_3$.
This we shall be able to do when,
critically, $p$ has the appropriate value. We will let
$p$ depend on $l$ and $k$, and the choice of $p$ is described
in more detail in Section \ref{sec-choicep}. Looking ahead, we
shall assume
    \beq
    \label{prange}
    k^{-3/2}\ll p\leq 10\frac{\ln k}{k}.
    \eeq
It will be convenient to subdivide the
possible $p$ into three regimes:

\ben
\item Very Large: $p\gg k^{-1}$ and $p\leq 10\frac{\ln k}{k}$
\item Large: $p=\frac{c}{k}$ for some $c>0$.
\item Small: $k^{-3/2}\ll p\ll k^{-1}$.
\een
In each of these cases, we will write $p=\frac{c}{k}$,
where $c\rightarrow 0$ when $p$ is small, and $c\rightarrow \infty$ when
$p$ is very large. The remainder of this section
is organized as follows. In Section \ref{sec-choicep}, we define how
to choose $p$ appropriately, and we show that the above
three regimes for $p$ correspond to the three regimes of $l$
given earlier. In Section \ref{twowalks}, we investigate
two walk problems, and relate the probability of $\TREE$ to the
probability that these two walks do not revisit their starting point 0.
In Section \ref{sec-gaussian}, we show that both $M$, and, more importantly,
$M^*$ obey a central limit theorem. Finally, in Section \ref{sec-asyC},
we state the consequences of our results concerning $\P[\TREE]$ and the
asymptotic normality of $M^*$ for $C(k,l)$.

\subsection{The Choice of Tilt}
\label{sec-choicep}
Let $\mu,\sig^2$ denote the mean
and variance of $M$.  Both of these have closed forms as a
function of $p$.
We have the exact calculation
    \beq\label{exactmu}
    \mu = (k-1)[\frac{k-1}{2}-E[T_1]] = (k-1)\Big[\frac{k}{2}-
    \frac{1-(k+1)p(1-p)^k-(1-p)^{k+1}}{p(1-(1-p)^k)}\Big].
    \eeq
We choose $p$ to satisfy the equation
    \beq
    \label{setp} p\mu = l.
    \eeq
We can show from Calculus that $\mu=\mu(p)$ is an increasing
function of $p$ and so (\ref{setp}) will have a unique solution.
The asymptotics depends on the regime. If $p=\frac{c}{k}$, then
    \beq\label{asymmularge}
    \boxed{p\mu \sim f_1(c)k \qquad
    \mbox{ and }
    \qquad \sig^2 \sim f_2(c)k^3.}
    \eeq
with
    \eq f_1(c) =
    c\left[\frac{1}{2}-\frac{1-(c+1)e^{-c}}{c(1-e^{-c})}\right] .
    \en
and, setting $\kappa=c(1-e^{-c})^{-1}$,
    \eq
    f_2(c) = \kappa\left[e^{-c}[-c^{-1}-2c^{-2}-2c^{-3}]+2c^{-3}
\right] -\left(\kappa[e^{-c}(-c^{-1}-c^{-2})+c^{-2}]\right)^2.
\en
In particular, for $p$ small, and using that $f_1(c)\sim \frac{c^2}{12}$
and $f_2(c)\sim \frac{1}{12}$ when $c\downarrow 0$,
    \beq\label{asymmu}
    p\mu = \frac{k^3p^2}{12} + O(k^4p^3)
    \qquad
    \mbox{ and }
    \qquad
    \sig^2\sim \frac{k^3}{12},
    \eeq
while for $p$ very large, and using that $f_1(c)\sim \frac{c}{2}$
and $f_2(c)\sim c^{-2}$ when $c\uparrow \infty$,
    \beq\label{asymmuverylarge} p\mu \sim  \frac{pk^2}{2}
    \qquad
    \mbox{ and }
    \qquad
    \sig^2 \sim \frac{k}{p^2}.
    \eeq
We see that the three regimes of $p$ do indeed correspond to the
three regimes of $l$, as we show now.
Indeed, for $l$ small, we have that $p\mu \sim \frac{k^3p^2}{12}=l$
is equivalent to
    \beq\label{asympsmall}
    p \sim k^{-3/2}\sqrt{12l}.
    \eeq
On the other hand, for $l$ large, if $l\sim k\beta$, then
    \beq\label{asymplarge}
    p\sim \frac{c}{k}
    \qquad
    \mbox{ with }
    \qquad
    f_1(c)=\beta,
    \eeq
while for $l$ very large,
    \beq\label{asympverylarge}
    p\sim 2lk^{-2}.
    \eeq
The asymptotics in (\ref{asymmularge}) with $p\sim \frac{c}{k}$ can be
found by approximating $k^{-1}T_j$ by the continuous truncated exponential
distribution over $[0,1]$, which has density $ce^{-cx}/(1-e^{-c})$.

\subsection{Two Walks}
\label{twowalks}
We define two basic walks.  In application the $Z_i,Z_i^{\sR}$ below will be
random variables of various sorts and so $\ESC$, $\ESC_{\sL}$, $\ESC^{\sR}$,
$\ESC_{\sL}^{\sR}$ become events.

\begin{definition}\label{leftwalk} Let $Z_1,Z_2,\ldots$ be nonnegative integers.
The {\em leftwalk} is defined by the initial condition $Y_0=1$ and the recursion
$Y_i=Y_{i-1}+Z_i-1$ for $i\geq 1$.  The escape event, denoted $\ESC$, is
that $Y_i>0$ for all $i\geq 1$.  The event $\ESC_{\sL}$, or escape until
time $L$, is that $Y_i>0$ for $1\leq i\leq L$.
\end{definition}

It shall often be convenient to count the bins ``from the right."  Let
$Z_i^{\sR}$, $1\leq i\leq k$, denote the number of balls in the $k-i+1$-st
bin.  Set $Y_0^{\sR}=0$ and $Y_i^{\sR}=Y_{i-1}^{\sR}+1-Z_i^{\sR}$ so that
$Y_i^{\sR}=Y_{k-i+1}$.  Then $\TREE$ becomes
    \beq
    \label{backtreeY}
    Y_t^{\sR} > 0 \mbox{ for } 1\leq t\leq k-1,
    \eeq
or, alternatively,
    \beq
    \label{backtree}
    Z_1^{\sR}+\ldots+Z_t^{\sR} \leq t-1 \mbox{ for } 1\leq t\leq k-1.
    \eeq
We shall generally use the superscript $R$ when examining bins from the
right.  In particular, we set $i^{\sR}=k-i+1$ so that bin $i^{\sR}$ is the
$i^{\rm th}$ bin from the right.

\begin{definition}\label{rightwalk}  Let $Z_1^{\sR},Z_2^{\sR},\ldots$ be nonnegative
integers.  The {\em rightwalk} is defined by the initial condition $Y_0^{\sR}=0$
and the recursion
$Y_i^{\sR}=Y_{i-1}^{\sR}+1-Z_i^{\sR}$ for $i\geq 1$.  The escape event, denoted $\ESC^{\sR}$, is
that $Y_i^{\sR}>0$ for all $i\geq 1$.  The event $\ESC_{\sL}^{\sR}$, or escape until
time $L$, is that $Y_i^{\sR}>0$ for $1\leq i\leq L$.
\end{definition}

We allow $Z_i,Z_i^{\sR}$ to be defined only for $1\leq i\leq L$
in which case $Y_i,Y_i^{\sR}$ are
defined for $0\leq i\leq L$ and $\ESC_{\sL},\ESC_{\sL}^{\sR}$
are well defined.
Indeed, our main results will be for these walks of length $L$, the
infinite walks shall be a convenient auxiliary tool.

When $k-1$ balls are placed into $k$ bins with tilt $p$ and
$Z_i$ is the number of balls in bin $i$ the event $\ESC_{\sL}$ is that
$Y_t>0$ for $1\leq t\leq L$.  Letting $Z_i^{\sR}$ be the number of balls
in bin $i^{\sR}$ the event $\ESC_{\sL}^{\sR}$ is that $Y_t^{\sR}>0$ for $1\leq t\leq L$.

\begin{definition} \label{leftrightmiddle}
Given $L<\frac{1}{2}k$, we call bins with $1\leq i\leq L$
the left side; bins with $1\leq i^{\sR} \leq L$ the right side; and all other
bins the middle.  \end{definition}

Now consider the tilted balls into bins formulation of Section \ref{statements}.
Set
    \beq\label{setlams}
    \lam = (k-1)\frac{p}{1-(1-p)^k}\qquad\mbox{ and }
    \qquad
    \lam^{\sR}=(k-1)
    \frac{p(1-p)^{k-1}}{1-(1-p)^k},
    \eeq
so that $\lam,\lam^{\sR}$ are the expected number of balls in the
leftmost and rightmost bin respectively.
When $p=\frac{c}{k}$, the asymptotics of $\lam$ and $\lam^{\sR}$
are given by
    \eq
    \lbeq{lambdadef}
    \boxed{\lam \sim \frac{c}{1-e^{-c}}\qquad
    \mbox{and}\qquad \lam^{\sR}\sim \frac{ce^{-c}}{1-e^{-c}}.}
    \en
In particular, for $p$ very large, $\lam \ra \infty$ and $\lam^{\sR}\sim 0$,
while for $p$ small, $\lam = 1+ \frac{pk}{2}(1+o(1))$ and
$\lam^{\sR}= 1- \frac{pk}{2}(1+o(1))$.

\begin{theorem}\label{probtreethm}
Let $\ESC$ be given by Definition \ref{leftwalk}, with all $Z_i\sim \Po(\lam).$
Let $\ESC^{\sR}$ be given by Definition \ref{rightwalk}, with all $Z_i^{\sR}\sim
\Po(\lam^{\sR})$.  Let $p$ be in the range given
by (\ref{prange}).  Then
    \beq\label{probtree}
    \Pr[\TREE] \sim \Pr[\ESC]\Pr[\ESC^{\sR}].
    \eeq
\end{theorem}

We may naturally interpret Theorem \ref{probtreethm} as saying that
the event $\TREE$ is  asymptotically equal to the probability that
the left and right sides satisfy the conditions imposed by $\TREE$.
For $i$ small, $Y_i$ behaves like a leftwalk with
$Z_i$ being roughly $\Po(\lam)$ and $Y_i^{\sR}$ behaves like a rightwalk
with $Z_i^{\sR}$ being roughly $\Po(\lam^{\sR})$. The proof
of Theorem \ref{probtreethm} is deferred to Section \ref{section2}.

The left and right walks with $Z_i,Z_i^{\sR}$ independent Poisson have been well
studied.
Let $Z\sim \Po(1+\eps)$.  Let $Z_i\sim Z$ all $i$, independent.
Consider the leftwalk as given by Definition \ref{leftwalk}.

\begin{theorem}\label{t1}
$\Pr[\ESC]=y$ where $y$ is the unique real number in $(0,1)$
such that $e^{-(1+\eps)y} = 1-y$.  Further, if $y=y(\eps)$, then
$y\leq 2\eps$ for all positive $\eps$ and $y\sim 2\eps$
as $\eps\ra 0^+$.
\end{theorem}

\proof
We use that there is a bijection between random walks with i.i.d.\ steps
with distribution $\Po(\lambda)-1$ and Galton-Watson trees with
offspring distribution $\Po(\lambda)$. This bijection is such
that random walks that never return to the origin are mapped
to branching process configurations where the tree is infinite.
For the latter, we have that the probability is the survival probability
of the branching process. The extinction probability $x$ satisfies
    \eq
    e^{-\lambda (x-1)} = x.
    \en
Therefore, for the survival probability $y=1-x$, we obtain
    \eq
    \lbeq{yimpeq}
    e^{-(1+\eps)y} = 1-y.
    \en
The inequality $y(\eps)\leq 2\eps$ and the asymptotics $y(\eps)\sim 2\eps$
are elementary calculus exercises.
\qed
\vskip0.5cm

\noindent
Next let $Z^{\sR}\sim \Po(1-\eps)$.
Let $Z_i^{\sR}\sim Z^{\sR}$ all $i$, independent.
Consider the rightwalk as given by Definition \ref{rightwalk}.
Then we can identify the probability of $\ESC^{\sR}$ exactly as follows:

\begin{theorem}\label{t2}
$\Pr[\ESC^{\sR}]=\eps$.
\end{theorem}

\proof
Consider an infinite walk starting at zero
with step size $1-P$ where $P$ is Poisson with mean $1-\eps$.
Here, $\eps\in (0,1)$ but we do {\em not} assume $\eps\ra 0$.
We claim $\Pr[\ESC^{\sR}] = \eps$ precisely. Take an infinite random walk,
$0=Y_0,Y_1,Y_2,\ldots$ and let $W_n$ be the number of
$i$,  where $0\leq i<n$,  for which $S_t=Y_{i+t}-Y_i$ for $t\geq 0$
never returns to zero, i.e., the number of
$i$, $0\leq i<n$ for which $Y_j>Y_i$ for all $j>i$.

For each $i$ this has probability $\ah$ of occurring
so  that by linearity of expectation $\E[W_n]=n\ah$.
Let $V_n$ be the minimum of $Y_j$ for $j\geq n$.  Then, by definition
$W_n=\max[V_n,0]$. Indeed, for each $0\leq j< V_n$, let
$i=i(j)$ be the maximal $i$, $0\leq i < n$ for which $Y_i=j$.  These
are precisely the $i$ for which the walk beginning at time $i$ has
the desired property.  Thus $n\ah = \E[\max[V_n,0]]$.  So far everything
is exact and now it follows from the fact that the random walk
has positive drift that
    \eq
    \lim_{n\ra\infty} \frac{\E[\max[V_n,0]]}{n} = \eps.
    \en
\qed
\vskip0.5cm

\noindent
Suppose $p=\frac{c}{n}$ with $c>0$ fixed.  Then $\Pr[\ESC]=y$ where
$e^{-\lam y}=1-y$ by \refeq{yimpeq} and $\lam$ is given by \refeq{lambdadef}
so that $\Pr[\ESC] \sim 1-e^{-c}$.  Theorem \ref{t2} gives  $\Pr[\ESC^{\sR}] =1-\lam^{\sR}$
where $\lam^{\sR}$ is given by \refeq{lambdadef} so that
$\Pr[\ESC^{\sR}] \sim \frac{1-(c+1)e^{-c}}{1-e^{-c}}$.
The asymptotics of Theorem \ref{probtreethm} are then that for $p=\frac{c}{k}$,
    \beq
    \label{probtreelarge}
    \boxed{\Pr[\TREE] \sim 1-(c+1)e^{-c}.}
    \eeq
In particular, for $p$ small
    \beq\label{probtreesmall} \Pr[\TREE] \sim \frac{p^2k^2}{2},
    \eeq
while for $p$ very large
    \beq\label{probtreeverylarge}
    \Pr[\TREE] \sim 1.
    \eeq

\subsection{The Limiting Gaussian}
\label{sec-gaussian}
In this section, we give an asymptotic normal law for $M^*$ and the
consequent asymptotics of $A_3$.
For $M$, by the fact that the $T_j$ are independent,
Esseen's Inequality gives that $M$
is asymptotically Gaussian with mean $\mu$ given in
(\ref{exactmu}) and variance $\sigma^2$. Therefore, for any fixed real $u$

    \beq\label{MGaussian}
    \Pr[M \leq \mu + u\sig ]\sim \int_{-\infty}^u \frac{1}{\sqrt{2\pi}}
    e^{-t^2/2}dt.
    \eeq

\begin{theorem}\label{main2}  Let $M^*$ be given by Definition \ref{def-M}.
Then for any fixed real $u$
    \beq\label{MGaussian1}
    \Pr[M^* \leq \mu + u\sig]\sim \int_{-\infty}^u \frac{1}{\sqrt{2\pi}}
    e^{-t^2/2}dt.
    \eeq
\end{theorem}

Here, importantly, $\mu$ is given by (\ref{exactmu}), the expectation
of the unconditioned $M$.  Theorem \ref{main2} then has the natural
interpretation that conditioning on $\TREE$ does not change the
asymptotic distribution of $M$.  The proof of Theorem \ref{main2} is
deferred to Section \ref{section2}.

We next use Theorem \ref{main2} to determine the asymptotics of $A_3$.
For this, we define $\sig_{\sY}$ by
    \beq
    \label{defsigY}
    \sig_{\sY}^2 = p\mu + p^2\sig^2.
    \eeq
\begin{prop}\label{Mpislthm} With $p$ given by (\ref{setp}) and
$\sig_{\sY}$ given by (\ref{defsigY}), whenever $p^2\sig^2=O(p\mu)$,
    \beq\label{Mpisl}
    \Pr[\BIN[M^*,p]=l]
    \sim \sig_{\sY}^{-1}(2\pi)^{-1/2}.
    \eeq
\end{prop}

\proof We require only the asymptotic Gaussian distribution
(\ref{MGaussian}). Using infinitesimals, the probability that $\sig^{-1}(M^*-\mu)\in [u,u+du]$
is $(2\pi)^{-1/2}e^{-u^2/2}du$.  Given that $M^*$ is in this range and using
that $l=p\mu$, we have
    \eq\Pr[\BIN[M^*,p]=l]\sim \Pr[\BIN[\mu+u\sig,p]=p\mu].
    \en
Note that the mean of $\BIN[\mu+u\sig,p]$ is equal to
$p(\mu+u\sig)$ and its variance is
    \eq
    p(1-p)(\mu+u\sig)\sim p(\mu+u\sig)\sim p\mu=l,
    \en
where in the latter equality, we use that $p\sig=O(\sqrt{p\mu})=O(\sqrt{l})=o(l)=o(p\mu)$.
Therefore, $p\mu$ is $\sim u(p\sig)(p\mu)^{-1/2}$ standard deviations off the mean $p\mu +
pu\sig$ and so this probability is $\sim (2\pi l)^{-1/2}\exp[-\frac{u^2}{2}
\frac{\sig^2p^2}{l}]$. Since $p^2\sig^2=O(p\mu)$, the values as
$u\ra \pm\infty$ are negligible and
    \beq\label{a12}
    \Pr[\BIN[M^*,p]=l] \sim (2\pi l)^{-1/2}\int_{-\infty}^{+\infty}
    e^{-\frac{u^2}{2}\frac{\sig_{\sY}^2}{l}}du,
    \eeq
which gives (\ref{Mpisl}).
\qed
\vskip0.5cm

\noindent
Again we can look at the asymptotics (we won't need finer expressions)
in the different regimes.
\ben
\item When $l$ is small,  then  $p\mu=l\sim p^2\sig^2$ and
    \beq\label{sigYsmall}
    \sig_{\sY}^2\sim 2l.
    \eeq
\item When $l$ is large, say $l\sim k\beta$,  then  $\sig^2p^2\sim c^2f_2(c)k$
with $c$ satisfying $f_1(c)=\beta$ as in (\ref{asymplarge}) so that
    \beq
    \label{sigYlarge}
    \sig_{\sY}^2 \sim l(1+c^2f_2(c)\beta^{-1}).
    \eeq
\item When $l$ is very large,  then  $\sig^2p^2=o(p\mu)$ and
    \beq
    \label{sigYverylarge}
    \sig_{\sY}^2 \sim l.
    \eeq
\een

\section{Asymptotics for $C(k,l)$}
\label{sec-asyC}
We now use the results in the previous section, in particular
Theorems \ref{exactthm}, \ref{probtreethm} and \ref{main2},
to derive the asymptotics for $C(k,l)$.
Indeed, for $p$ given by (\ref{setp}), the asymptotics of all terms in
(\ref{x}) are known except $C(k,l)$.  Hence we can solve for
the asymptotics of $C(k,l)$.
While Theorem \ref{exactthm} and the auxiliary results allow us
to find the asymptotics of $C(k,l)$ in theory, some of the technical
work can be challenging.  Here we indicate some of the major cases.

It shall be helpful not to use the precise $p$ given by (\ref{setp}).
Recall that Theorem \ref{exactthm} holds for any value of $p$.  For
the moment let us write $\mu=\mu(p)$, $\sig=\sig(p)$, $A_2=A_2(p)$ and
$A_3=A_3(p)$ to emphasize this dependence.

\begin{lemma}\label{easyapprox} Let $p_0$ be the value of $p$
satisfying (\ref{setp}), i.e., $p_0\mu(p_0)=l$.  Let $p$ be such
that $p\sim p_0$ and $p\mu(p) = l+ o(l^{1/2})$.  Then $A_2(p)\sim
A_2(p_0)$, $\sig(p)\sim \sig(p_0)$ and $A_3(p)\sim A_3(p_0)$.
\end{lemma}

\proof
The asymptotics of $A_2(p)=\Pr[\TREE]$ are given by (\ref{probtreelarge})--
(\ref{probtreeverylarge}) and clearly have this
property.  Similarly the formulae (\ref{asymmularge}),(\ref{asymmu})--
(\ref{asymmuverylarge}) show that $\sig(p)\sim \sig(p_0)$.  An
examination of the proof of Proposition \ref{Mpislthm} gives that the
asymptotics (\ref{Mpisl}) of  $\Pr[\BIN[M^*,p_0]=l]$ apply to
$\Pr[\BIN[M^*,p]=l]$ as long as $l-p\mu(p) = o(l^{1/2})$ as the integral
(\ref{a12}) remains the same.
\qed

\subsection{$l$ Small}
    \begin{theorem} When $l=o(k^{1/2})$,
    \eq
    \label{hooray1}
    C(k,l) \sim k^{k-2}k^{3l/2}(e/12l)^{l/2}(3\pi^{-1/2})l^{1/2}.
    \en
    \end{theorem}

\proof Setting $p=k^{-3/2}\sqrt{12l}$ we have from (\ref{asymmu}) that
$p\mu = l+O(k^4p^3)$ and $O(k^4p^3)=O(k^{-1/2}l^{3/2})=o(l^{1/2})$
as $l=o(k^{1/2})$.  Lemma \ref{easyapprox} then allows us to use this
$p$ with $A_2,A_3$ given by the $p$ of (\ref{setp}). We start with
the exact formula
    \eq
    \label{Crew}
    C(k,l)=A_1A_2A_3 p^{-(k+l-1)}(1-p)^{-{k\choose 2}+(k+l-1)}.
    \en
By (\ref{probtreesmall}), we have
    \eq
    A_2=\Pr[\TREE]\sim \frac{1}{2}(kp)^2=6lk^{-1}.
    \en
We further have $p^2\sig^2+p\mu \sim 2l\sim \frac{1}{6}\lam^2$
so that Proposition \ref{Mpislthm} gives
    \eq
    A_3=\Pr[\BIN[M^*,p]=l]\sim (2\pi)^{-1/2}(2l)^{-1/2}.
    \en
We have to be quite careful with the asymptotics of the exact formula
$A_1=[1-(1-p)^k]^{k-1}$ of Theorem \ref{exactthm}.  We have
    \beq
    \label{techa11}
    1-(1-p)^k = pk(1-\frac{pk}{2}+\frac{p^2k^2}{6}-\frac{p^3k^3}{24} + O(p^4k^4)),
    \eeq
and
    \beq
    \label{techa12}
    \ln(1-\frac{pk}{2}+\frac{p^2k^2}{6}-\frac{p^3k^3}{24} + O(p^4k^4))=
    -\frac{pk}{2}+\frac{p^2k^2}{6}-\frac{p^3k^3}{24}-\frac{p^2k^2}{8}+
    \frac{p^3k^3}{12}-\frac{p^3k^3}{24}+O(p^4k^4)
    =-\frac{pk}{2} +\frac{l}{2k}+o(k^{-1}),
    \eeq
as long as $\lam=o(k^{1/4})$ which occurs if $l=o(k^{1/2})$.
This gives
    \beq
    \label{techa13}
    A_1\sim p^{k-1}k^{k-1}e^{-pk^2/2}e^{l/2}.
    \eeq
We further have the asymptotics
    \beq
    \label{techa14}
    p^{k+l-1}\sim p^{k-1}(12l)^{l/2}k^{-3l/2},
    \eeq
and
    \beq
    \label{techa15}
    (1-p)^{k^2/2} \sim e^{-k^2p/2},
    \eeq
while
    \eq
    (1-p)^{k+l-1} \sim 1.
    \en
Then Theorem \ref{exactthm} puts everything together and yields (\ref{hooray1}).
\qed

\subsection{$l$ Large}
In this section, we take $l$ such that $\beta\equiv l/k$ is uniformly bounded and
uniformly positive, and investigate the scaling of $C(k,l)$ in this range. We state
the result uniformly in $\beta\in [\eps k, \eps^{-1} k]$, since we cannot fix $\beta$
due to the fact that $l=\beta k$ need to be an integer. The main result
is as follows:

\begin{theorem}
\label{thm-Cllarge}
Let $\eps>0$ be arbitrary and fixed.  Then as
$k\ra\infty$ and $l=l(k)\in [\eps k,\eps^{-1} k]$,
    \beq
    C(k,k\beta)\sim A \cdot B^k \cdot k^{(1 + \beta)k} \cdot k^{-3/2},
    \eeq
where $\beta=l/k$, $c$ is the solution to
\beq \label{csol}
    e^{-c} = {\frac{2(\beta + 1)-c}{2(\beta +1) + c}},
\eeq
    \beq
    A = {\frac{c(c - 2 \beta )}{\sqrt{8\pi\beta(1+c^2 f_2(c)/\beta)}}}
    e^{-c(\beta/2 + 1)},
    \eeq
and
    \beq
    B = {\frac{2}{c^\beta \sqrt{4(\beta+1)^2 - c^2}}}.
    \eeq
\end{theorem}

\proof Let $l$ satisfy that $l=l(k)\in [\eps k,\eps^{-1} k]$.
Then the $p$ of (\ref{setp}) satisfies
    \beq p= \frac{c}{k} + O(k^{-2}),
    \eeq
where c is the solution to (\ref{csol}).
Changing $p$ by an additive $O(k^{-2})$ term changes $p\mu(p)$ by
$O(1)$.  Lemma \ref{easyapprox} allows us to set $p=\frac{c}{k}$
with $A_2,A_3$ the same as for that $p$ given by (\ref{setp}).
We get $C(k,\beta k)$ from the equation
    \beq
    A_1 A_2 A_3 =  C(k,\beta k) p^{k+\beta k -1} (1-p)^{{k \choose 2}-k-\beta k +1}
    \eeq
Here, taking care to note  that the asymptotics $(1-p)^k\sim e^{-c}$ are {\em not}
sufficiently precise to give the asymptotics of $A_1$, we find
   \beq  A_1 \sim (1 - e^{-c}(1-\frac{c^2}{2k}))^{k-1}
    \sim (\frac{2c}{2(\beta +1) + c})^{k-1}
    e^{((\beta +1)/c + 1/2)c^2/2},
    \eeq
while
    \beq
    A_2 \sim 1 -(c+1)e^{-c}, \qquad \text{and}\qquad
    A_3 = \frac{1}{\sqrt{2\pi k\beta (1+c^2 f_2(c)/\beta)}}.
    \eeq
Further,
    \beq
    p^{k(1+\beta) - 1} \sim (\frac{c}{k})^{k(1+\beta) - 1}
    ,
    \eeq
while
    \beq (1-p)^{{k \choose 2}-k(1+\beta) + 1}
    \sim e^{-kc/2-c_1/2-c^2/4+c(\beta +1.5)},
    \eeq
Solving and employing uniformity of convergence we obtain Theorem \ref{thm-Cllarge}.
\qed

%
%

\subsection{$l$ Very Large}

As a third example suppose $l=\lfloor ck\ln k\rfloor$. We prove the following result:
    \begin{theorem}
    For $l=\lfloor ck\ln k\rfloor$ with $c>\frac 12$,
        \eq
        C(k,l) \sim {{\frac{k(k-1)}{2}}\choose {k+l-1}}.
        \en
    \end{theorem}

This has the interpretation that the proportion of graphs with $k$ vertices
and $k+l-1$ edges which are connected is asymptotically one, or that the
probability that a random graph with $k$ vertices and $k+l-1$ edges is
connected is asymptotically one.  As such, this
is immediate from a classic results of Erd\H{o}s and R\'enyi \cite{ER}.

\proof
We again start from (\ref{Crew}). Then
(\ref{asympverylarge}) gives that $p\sim 2lk^{-2}$,
which implies $A_1\sim A_2\sim 1$.
(Note that $A_1\sim 1$ fails for $c<\frac 12$.)
Further, Proposition \ref{Mpislthm}
with the asymptotics in (\ref{sigYverylarge})
gives $A_3\sim (2\pi l)^{-1/2}$. It shall be convenient to rewrite this
as $A_3\sim (2\pi (k+l-1))^{-1/2}$. We conclude that
    \eq
    C(k,l)\sim (2\pi (k+l-1))^{-1/2} p^{-(k+l-1)}(1-p)^{-{k\choose 2}+(k+l-1)}
    = (2\pi B)^{-1/2} p^{-B}(1-p)^{-A+B},
    \en
where we abbreviate $A={k\choose 2}, B=k+l-1$.
However, this is not a sufficiently precise approximation of $p$ to give
the asymptotics of $C(k,l)$.  Rather, in the region $p\sim 2lk^{-2}$, the exact
expression (\ref{exactmu}) can be rewritten as follows:

    \begin{lemma}
    \label{lem-approxfullT}
        \eq
        \mu = {k\choose 2} - (k-1)p^{-1} -\frac{k(k-1)(1-p)^k}{1-(1-p)^k}.
        \en
    \end{lemma}

\proof This is a simple calculation.
\qed
\vskip0.5cm
\noindent
By Lemma \ref{lem-approxfullT}, and using that $1-(1-p)^k\sim 1$, we obtain that
    \eq
    \mu = {k\choose 2} - (k-1)p^{-1} + O(k^2(1-p)^k).
    \en
As we have required
from (\ref{setp}) that $p\mu =l$, we have
    \eq
    \lbeq{pmuasy}
    p{k\choose 2} = k+l-1 +O(pk^2(1-p)^k).
    \en
We note that
    \beq\label{aae}
    {A\choose B}p^B(1-p)^{A-B}=\Pr[\BIN[A,p]=B]
    \sim(2\pi B)^{-1/2},
    \eeq
where the latter equality holds by the local central limit theorem for the
binomial distribution whenever $B-pA=o(\sqrt{p(1-p)A})$.
Note that by \refeq{pmuasy},  with $A={k\choose 2}$ and $B=k+l-1$,
    \eq
    B-pA= (k+l-1)-p{k\choose 2}=O(pk^2(1-p)^k)=o(\sqrt{p(1-p)A})=o(k \sqrt{p}),
    \en
precisely when
    \eq
    (1-p)^k \sqrt{pk^2} =o(1).
    \en
Since $p\sim 2lk^{-2}$, this holds precisely when
    \eq
    \sqrt{l}e^{-\frac{2l}{k}}=o(1),
    \en
which is true whenever $l=c\frac{\log{k}}{k}$ with $c>\frac 14$.
Therefore, in this case, Theorem \ref{exactthm} gives
    \beq
    \label{aaf}
    C(k,l)p^B(1-p)^{A-B}
    \sim (2\pi B)^{-1/2},
    \eeq
and thus we deduce
    \beq\label{aag}
    C(k,l) \sim {A\choose B}={{\frac{k(k-1)}{2}}\choose {k+l-1}}.
    \eeq
\qed
\vskip0.5cm

\noindent
We note that in principle it is possible to extend the above asymptotics to
other $l$ for which $\frac{l}{k}\rightarrow \infty$, using Lemma \ref{lem-approxfullT}
and more precise local central limit theorems for $\Pr[\BIN[A,p]=B]$.

\section{The Technical Theorems}\label{section2}
In this section,  we prove Theorems \ref{probtreethm} and \ref{main2}.  The
values $\lam,\lam^{\sR}$, the expected number of balls in the first
and last bins respectively, are given by (\ref{setlams}).
We start in Section \ref{easycases} with the easy
case where $p$ is large and very large. The remaining Sections
\ref{sec-hard1}--\ref{sec-hard9} are devoted to the hard case
where $p\gg k^{-3/2}$.

\subsection{The Easy Case: $p$ Very Large and $p$ Large}
\label{easycases}
We note that the arguments for the ``hard case" apply to the cases
where $p$ is large and very large as well. However, many of the
subtleties of the hard case can be avoided when
$p= \Omega(k^{-1})$. Here we give, without full details,
a simpler argument that works in these important cases.

First suppose $p\gg k^{-1}$.  Let $\FAIL_t$ be the event $Y_t\leq 0$.
For example, $\FAIL_1$ is the
event $Z_1=0$ which has probability $e^{-(1+o(1))\lam}$
which approaches zero.  The event $\FAIL_k$ is the
event $Z_1^{\sR}>0$ so $\Pr[\FAIL_k]\leq E[Z_1^{\sR}] = \lam^{\sR}\ra 0$.
In general, as each ball is dropped
independently $Z_1+\ldots+Z_t$ has distribution
$\BIN[k-1,\ah]$ where $\ah = \Pr[T_j\leq i]$ as given
by the distribution (\ref{dist}).  (Near the right side
it is easier to work with $\Pr[Y_t^{\sR}\leq 0]$.)
Chernoff bounds give that $\sum_{t=1}^{k}\Pr[\FAIL_t]\ra 0$ and
so $\Pr[\TREE]\ra 1$, giving Theorem \ref{probtreethm}.
Conditioning on an event that holds with probability
$1-o(1)$ cannot effect an  asymptotic Gaussian distribution
and so Theorem \ref{main2} follows immediately for the very large case.

 We next proceed with the case where $p$ is large.
Set $p=\frac{c}{k}$.  Note $\lam,\lam^{\sR}$ are given by (\ref{setlams}).
We split bins into left, right and middle by Definition \ref{leftrightmiddle}.
We set $L=\lfloor\ln^2k\rfloor$ for definiteness, though a fairly wide
range of $L$ would do. With $\FAIL_t$ as above, Chernoff bounds give
$\sum_{L<t\leq k-L}\Pr[\FAIL_t] = o(1)$.  With probability $1-o(1)$,
no bin on the left nor right side has more than $\ln^2k$ balls so
that the total number of balls on the left and right side is less
than $\ln^4k$. Thus, $\Pr[\TREE]$ is within $o(1)$ of the probability
that both sides have less than $\ln^4k$ balls and that the leftwalk
satisfies $\ESC_{\sL}$ and that the rightwalk
satisfies $\ESC_{\sL}^{\sR}$.

Let $Z_i^*\sim \Po(\lam)$, $1\leq i\leq L$, be independent.
Let $Z_i^{\sss R*}\sim \Po(\lam^{\sR})$, $1\leq i\leq L$, be independent.
Placing balls into the left and right
sides with these distributions with probability $1-o(1)$ both
left and right sides get less than $\ln^4k$ balls. Both
$\Pr[\ESC_{\sL}],\Pr[\ESC_{\sL}^{\sR}]$ are within $o(1)$
of $\Pr[\ESC],\Pr[\ESC^{\sR}]$ for the infinite walks and, as they are
now independent, $\Pr[\ESC_{\sL}\wedge \ESC_{\sL}^{\sR}]$ would be within
$o(1)$
of $\Pr[\ESC]\Pr[\ESC^{\sR}]$.
For any fixed nonnegative integers $x_1,\ldots,x_{\sL};x_1^{\sR},\ldots,x_{\sL}^{\sR}$
the probability that $Z_i=x_i$, $1\leq i\leq L$ and $Z_i^{\sR}=x_i^{\sR}$,
$1\leq i\leq L$ approaches the same
probability with the $Z_i,Z_i^{\sR}$ replaced by the independent
Poissons $Z_i^*,Z_i^{\sR *}$.  Hence, $\Pr[\TREE]$ is within $o(1)$
of $\Pr[\ESC]\Pr[\ESC^{\sR}]$, giving Theorem \ref{probtreethm}.

We next proceed with the central limit theorem Theorem \ref{main2}.
The proof of this result is more subtle, and we need to show that
both mean and variance are not affected by the conditioning.
Consider any fixed nonnegative integers $x_1,\ldots,x_{\sL};x_1^{\sR},\ldots,x_{\sL}^{\sR}$ so
that, with $x_i$ balls in bin $i$ and $x_i^{\sR}$ balls in bin $i^{\sR}$ the events
$\ESC_{\sL}$ and $\ESC_{\sL}^{\sR}$ both hold.
Set $m_{\sL}=x_1+\ldots+x_{\sL}$, $m_{\sR}=x_1^{\sR}+\ldots+x_{\sL}^{\sR}$ and
further assume $m_{\sL}< \ln^4k$ and $m_{\sR}<\ln^4k$.  Let
$M^{**}$ be the distribution of ${k\choose 2}-\sum T_x$ where we assume
that all remaining balls are placed in the middle bins with the truncated
geometric distribution.  Thus, the law of $M^{**}$ is the law of
$M^*$ conditioned on $m_{\sL}< \ln^4k$ and $m_{\sR}<\ln^4k$.
Let $\mu^{**}=E[M^{**}]$. Then, the following proposition shows that
the conditioning does not affect the mean too much:

\begin{prop}\label{easyclaim} $\mu^{**}-\mu = O(k\ln^4k)$.
\end{prop}

\proof
Lets call these distributions fixededge and unrestricted respectively.
There are two differences between these distributions.  First, the
$m_{\sL}+m_{\sR}$ balls are explicitly placed in the fixededge
distribution.  The difference in expectation for any particular ball
can be at most $k$ so the total difference for these less than
$\ln^4k$ balls is less than $k\ln^4k$.  For the other balls the
distinction is between the truncated geometric and the unrestricted
distribution.  Let $Y_{\sss T},Y_{\sss U}$ be the placement of a single ball in
these two distributions.
Consider the experiment of selecting $Y_{\sss T}$
from the unrestricted distribution and then reassigning it with
the truncated geometric if it did not land in a middle bin.
With this linkage we have $Y_{\sss T}\neq Y_{\sss U}$ only when the reassignment
is made which occurs with probability $O(k^{-1}\ln^2k)$.  When it
does occur the values are, as always, within $k$.  Hence the difference
in the expectations is $O(\ln^2k)$.  The total  difference for all
(at most $k$) of these balls is then $O(k\ln^2k)$.  Thus
$\mu^{**}-\mu = O(k\ln^4k)+O(k\ln^2k)$ giving Proposition
\ref{easyclaim}.
\qed
\vskip0.5cm

Now we claim that $M^{**}$ satisfies the asymptotic Gaussian (\ref{MGaussian}).
We may write $M^{**}=\ah - \sum T_j^{**}$ where $\ah$ is a constant which
depends on the fixed placement, the sum ranges over those $j$ for which
ball $j$ goes into the middle, and $T_j^{**}$ has the distribution of
$T$ given by (\ref{dist}) conditioned on it being in the middle.  We
claim $M^{**}$ has variance $\sim \sig^2$ with $\sig^2$ given by
(\ref{asymmularge}).  For $M,M^{**}$ the variance comes from
the independent $T_j,T_j^{**}$ respectively.  There are $k-1$ and
$\sim k$ terms respectively.  The variance of each $T_j$ and each
$T_j^{**}$ is $\sim f_2(c)k^2$. An easy way to see this is that
$k^{-1}T_j^{**}$ has the asymptotic continuous distribution on
$[0,1]$ with density $e^{-cx}/(1-e^{-c})$,  which is the asymptotic
law of $k^{-1}T_j$ when $k\rightarrow \infty$.
From Esseen's Inequality, $M^{**}$ is asymptotically Gaussian with mean
$\mu^{**}$ and variance $\sim \sig^2$.  Since $\mu-\mu^{**}=
O(k\ln^4k)=o(k^{3/2})$,
$M^{**}$ is asymptotically Gaussian with the original $\mu,\sig^2$.

Finally, we consider $M^{*}$.  In the unconditioned placement of balls
the probability that either $m_{\sL}>\ln^4k$ or $m_{\sR}>\ln^4k$ was
$o(1)$.  We are now conditioning on $\TREE$ but we have already shown
that, in this regime, $\Pr[\TREE]$ is bounded away from zero.  Hence,
in the conditioned placement of balls the probability that either
$m_{\sL}>\ln^4k$ or $m_{\sR}>\ln^4k$ is still $o(1)$.
Therefore, excluding $o(1)$ probability, $M^{*}$ is a combination of
of distributions $M^{**}$, each of which is
asymptotically Gaussian with the same mean and variance.
Hence, $M^{*}$ is as well.  This completes the proof of Theorem
\ref{main2} in the case when $p$ is large.
\qed

\subsection{The Hard Case}
\label{sec-hard1}

In Section \ref{sec-hard1}--\ref{sec-hard9}, we study the general case where $p k^{3/2}
\rightarrow \infty$. Our arguments can be made considerably simpler when
$p$ is not too close to the lower bound $k^{-3/2}$.  When we present the
general results, we will indicate the
simplification when $p=k^{-1.4}$.
These simplifications actually work down
to $k^{-3/2}$ times a polylog factor.

We split the $k$ bins into left, middle and right sides as given
by Definition \ref{leftrightmiddle}.
We carefully choose $L$
so that
    \beq\label{Lbounds}
    (kp)^{-2}\ll L \ll k^{-1/2}p^{-1}
    \eeq
For example, when $p=k^{-1.4}$, we set $L=k^{0.85}$, far
away from both bounds of (\ref{Lbounds}).

Note that the lower bound of (\ref{prange}) on $p$ allows us to do this.
Also note that $k^{-1/2}p^{-1}\ll k$ so that
    \beq\label{Lbounds1}
    L\ll k.
    \eeq
We set
    \beq
    \label{seteps}
    \eps = \frac{pk}{2}=o(1).
    \eeq
since $p$ is small. A careful analysis of (\ref{dist}) gives that
    \beq\label{leftdist}
    \Pr[T_j=i] = \frac{1}{k}\left(1 + \eps+o(\eps)\right) \quad\mbox{ for } \quad 1\leq i\leq L,
    \eeq
and
    \beq\label{rightdist}
    \Pr[T_j=i^{\sR}] = \frac{1}{k}\left(1 - \eps+o(\eps)\right) \quad \mbox{ for }\quad
    1\leq i\leq L.
    \eeq
Roughly speaking, each bin on the left side will get $\Po(1+\eps)$
balls, while the bins on the right side will get $\Po(1-\eps)$ balls.
It shall turn out that the event $\TREE$ is dominated by the events of
(\ref{tree}) for $1\leq t\leq L$ and the events of (\ref{backtree})
for $1\leq t\leq L$.\\

\subsection{Scaling for Small Bias Walks}
\label{sec-hard3}

Mathematical physicists well understand that walks with a bias
$\eps=o(1)$ are naturally scaled by time $\eps^{-2}$.
Up to time $O(\eps^{-2})$ the walk behaves as if it had zero drift
and afterwards the drift takes over.  Propositions \ref{t3}--\ref{t4}
below investigate the probability of never returning to the starting point,
and are quite natural. We write $\Pr_{\eps}^*$ for the law where each
bin $1, 2, \ldots$ receives a Poisson($1+\eps$) number of balls.

\begin{prop}\label{t3}
If $\eps\ra 0^+$ and $L\ra \infty$ is such that $L\gg \eps^{-2}$,
then $\Pr_{\eps}^*[\ESC_{\sL}] \sim 2\eps$.
\end{prop}

\proof
As $\Pr_{\eps}^*[\ESC]\sim 2\eps$ it suffices to show $\Pr_{\eps}^*[\neg \ESC\wedge
\ESC_{\sL}] = o(\eps)$.

In the simpler case when $p=k^{-1.4}$ so $\eps=k^{-0.4}/2$ and
$L=k^{0.85}$, we can bound $\Pr_{\eps}^*[\ESC_{\sL}]$ by the sum over $t>L$
of $\Pr_{\eps}^*[Z_1+\ldots+Z_t < t]$.  Here $Z_1+\ldots+Z_t\sim \Po(t(1+\eps))$.
Basic Chernoff bounds show that this probability is so low and drops
so fast that summed over all $t>L$ it is $o(\eps)$.  Indeed, it is
of the form $\exp[-k^{c+o(1)}]$ for some positive constant $c$.
Now we extend the proof to the small $p'$s for which
$pk^{3/2}\rightarrow \infty$.

Consider the infinite walk and let $W$
be the number of $t\geq L$ such that $Y_t \leq \frac{L\eps}{2}$.
Parametrize $t=Lx$.  Then
    \eq
    \Pr~\!\!_{\eps}^*[ Y_t \leq \frac{L\eps}{2}] \leq
    \Pr~\!\!_{\eps}^*[\Po(Lx(1+\eps))\leq Lx+\frac{Lx\eps}{2}].
    \en
Basic Chernoff bounds give that this is at most
$\exp[-(Lx\eps)^2/8(Lx(1+\eps))] \leq \exp[-L\eps^2x/16]$.
Since $L\eps^2\gg 1$, this probability is $o(1)$ for every $x\in(0,1)$ fixed.
Therefore, by linearity of expectation, $E[W] = o(L)$.  Let $B$ be the event that
$Y_t=0$ for some $t\geq L$.  Then we  claim $E[W|B] \geq 0.98L$.

Indeed, consider the first such $t\geq L$ with $Y_{t}=0$.  Conditionally on
the history up to time $t$, we have $\Pr_{\eps}^*[Y_{t+s} \leq L\eps/2] \geq 0.99$ for
all $0\leq s\leq L(0.99)$.  As $E[W]\geq E[W|B]\Pr_{\eps}^*[B]$ we deduce that
$\Pr_{\eps}^*[B]=o(1)$.  Now $\ESC_{\sL}$ is an increasing event and $B$ is a decreasing
event, so that by the FKG inequality
    \eq
    \Pr~\!\!_{\eps}^*[\ESC_{\sL}\wedge B]\leq \Pr~\!\!_{\eps}^*[\ESC_{\sL}]\Pr~\!\!_{\eps}^*[B],
    \en

\noi so that
$\Pr_{\eps}^*[\ESC_{\sL}\wedge\neg \ESC] = \Pr_{\eps}^*[\ESC_{\sL} \wedge B] = o(\eps)$.
We conclude that
    \eq
    \Pr~\!\!_{\eps}^*[\ESC_{\sL}]=\Pr~\!\!_{\eps}^*[\ESC]+\Pr~\!\!_{\eps}^*[\ESC_{\sL}\wedge B]
    =\Pr~\!\!_{\eps}^*[\ESC]+ o\big(\Pr~\!\!_{\eps}^*[\ESC_{\sL}]\big),
    \en
so that $\Pr~\!\!_{\eps}^*[\ESC_{\sL}]\sim \Pr~\!\!_{\eps}^*[\ESC]$. Since, by
Theorem \ref{t1}, we have $\Pr~\!\!_{\eps}^*[\ESC]\sim 2\eps$, Proposition \ref{t3}
follows. \qed
\vskip0.5cm

\noindent
The next proposition gives a similar result for $\ESC_{\sL}^{\sR}$.
In its statement, we let $\Pr_{\sR,\eps}^*$ denote the probability law where each bin
$1, 2, \ldots, \infty$ receives a Poisson($1-\eps$) number of balls.

\begin{prop}\label{t4}
If $\eps\ra 0^+$ and $L\ra \infty$ is such that $L\gg \eps^{-2}$,
then $\Pr_{\sR,\eps}^*[\ESC_{\sL}^{\sR}] \sim \eps$.
\end{prop}

\proof Similar to the proof of Proposition \ref{t3}.
\qed
\vskip0.5cm

\noindent
We further require two small extensions:

\begin{cor}\label{t10}  Let $\eps\ra 0^+$ and $L\gg \eps^{-2}$.
Let $\lam_1,\ldots,\lam_{\sL}$ be such that all $\lam_i=1+\eps+o(\eps)$.
Let $Z_i\sim \Po(\lam_i)$ be independent and consider the leftwalk
defined by Definition \ref{leftwalk}.  Then $\Pr[\ESC_{\sL}]\sim 2\eps$.
\end{cor}

\proof For any fixed $\del>0$ we can sandwich this model between
one in which all $\lam_i=1+\eps(1-\del)$ and one in which all
$\lam_i=1-\eps(1-\del)$.  From Proposition \ref{t3} we bound
$\Pr[\ESC_{\sL}]$ between $\sim 2\eps(1-\del)$ and $\sim 2\eps(1+\del)$.
As $\del$ can be arbitrarily small this gives Corollary \ref{t10}.
\qed

\begin{cor}\label{t11}  Let $\eps\ra 0^+$ and $L\gg \eps^{-2}$.
Let $\lam_1^{\sR},\ldots,\lam_{\sL}^{\sR}$ be such that all $\lam_i^{\sR}=1-\eps+o(\eps)$.
Let $Z_i^{\sR}\sim \Po(\lam_i^{\sR})$ be independent and consider the rightwalk
defined by Definition \ref{rightwalk}.
Then $\Pr[\ESC_{\sL}^{\sR}]\sim \eps$.
\end{cor}

\proof Similar as the proof of Corollary \ref{t10}, now
using Proposition \ref{t4} instead of Proposition \ref{t3}.
\qed

\subsection{Poisson versus Fixed}  \label{versus}
It will often be convenient to start with a Poisson number of balls with
parameter $k$, rather than with precisely $k$ balls. Indeed, in the Poisson
case, the number of balls per bin are independent Poisson random variables,
which is often quite convenient in the analysis. Therefore, sometimes
we wish to compare the probability that an event holds when we use
a Poisson number of balls
to the probability that the event holds when we use a fixed number of balls.
In this section, we prove a result that allows us to compare these
probabilities, and which will in particular allow us to convert a probability
for the Poisson law into a statement for the probability of the event for a
fixed number of balls.

We first introduce some notation that allows us to make this comparison.
Consider an event $A$ that depends on a nonnegative integer variable $X$.
Let $g(\lam)$ be $\Pr[A]$ when $X$ has a Poisson distribution with mean
$\lam$.  Let $f(m)$ be $\Pr[A]$ when $X=m$.  (As an important example,
drop $X$ balls into bins $1,\ldots,L$ with left-tilt $p$.)  These are
related by the equality
    \beq\label{xy7}
    g(\lam)= \sum_{m=0}^{\infty} f(m)\Pr[\Po(\lam)=m].
    \eeq
Here we want to go from asymptotics of $g$ to asymptotics of $f$.
We would naturally want to say that $g(m)$ and $f(m)$ are quite
close. This is true when $f$ and $g$ are {\it increasing} or
{\it decreasing.}

We say $A$ is increasing if $f,g$ are
increasing; decreasing if $f,g$ are decreasing and monotone if one
of those hold.  For balls into bins models, an event $A$ is increasing when
$A$ keeps on holding when extra balls are added.  An event
$A$ is decreasing when $A^c$
is increasing.   In particular, $\ESC_{\sL}$, $\ESC_{\sL}^{\sR}$ are
increasing and decreasing respectively.  When $A$ is monotone and
$g$ is relatively smooth the following result allows us to derive the
asymptotics of $f$ from those of $g$:

\begin{lemma}\label{poissonfixed}
Let $\lam_1,\lam_2\ra \infty$ with $\lam_2=\lam_1+\omega \lam_1^{1/2}$
where $\omega\ra \infty$.  Suppose $g(\lam_1)\sim g(\lam_2)$.
Then: \\ If $A$ is increasing, then $f(\lam_1)\leq (1+o(1))g(\lam_2)$.
\\ If $A$ is decreasing, then $f(\lam_2)\leq (1+o(1))g(\lam_1)$.
\\ If $A$ is increasing, then $f(\lam_2)\geq (1+o(1))g(\lam_1)$.
\\ If $A$ is decreasing, then $f(\lam_1)\geq (1+o(1))g(\lam_2)$.
\end{lemma}

\proof
Assume $A$ is increasing.  Truncating (\ref{xy7}) to
$m\geq \lam_1$ gives $g(\lam_2)\geq f(\lam_1)\Pr[\Po(\lam_1)\geq \lam_2]$.
Chebyschev's Inequality gives that the probability is $1-o(1)$, giving
the first part of Lemma \ref{poissonfixed}.  Now we show the third
part.  Calculation gives that for $j\geq \lam_2$, $Pr[\Po(\lam_2)=j]\gg
Pr[\Po(\lam_1)=j]$.  Now consider the expansion (\ref{xy7}) for both
$\lam=\lam_1$ and $\lam=\lam_2$.  We bound
    \beq\label{wwe}
    \sum_{m\geq\lam_2} f(m)\Pr[\Po(\lam_1)=m] \ll
    \sum_{m\geq\lam_2} f(m)\Pr[\Po(\lam_2)=m] \leq g(\lam_2)\sim g(\lam_1),
    \eeq
so that
    \beq\label{wwi}
    g(\lam_1)\sim \sum_{m<\lam_2} f(m)\Pr[\Po(\lam_1)=m]  \leq f(\lam_2).
    \eeq
Statements two and four are similar.
\qed

In application we will deal with situations in which $g(\lam)$ is asymptotically
constant in an interval around $\lam_0$ of width $\gg \sqrt{\lam_0}$.  In
that case $f(m)\sim g(m)$ for all $m$ in that interval.

\subsection{The probability of $\TREE$ in the left, right and middle bins}
\label{sec-hard5}
In this section, we investigate the probabilities of $\TREE$ in the left,
right and middle bins. The main results are Propositions \ref{t24}, \ref{t25}
and \ref{middle}. In Sections \ref{sec-hard9}--\ref{sec-hard7},
these results, as well as Corollaries \ref{t10} and \ref{t11},
will be combined to prove Theorem \ref{probtreethm}.

We first use  Lemma \ref{poissonfixed} together with
the results in Corollaries \ref{t10} and \ref{t11} to
investigate the probabilities of $\ESC_{\sL}$ and of $\ESC_{\sL}^{\sR}$:

\begin{prop}\label{t24} Let $\eps\ra 0^+$ and $L\gg \eps^{-2}$.
Let $m=L(1+\eps+o(\eps))$.  Let $p_1,\ldots,p_{\sL}$ have all
$p_i = \frac{1}{L} + o(\frac{\eps}{L})$.  Let $f(m)$ denote the
probability of $\ESC_{\sL}$ when precisely $m$ balls are placed in
bins $1,\ldots,L$ according to this distribution.  Then $f(m)\sim 2\eps$.
\end{prop}

\proof
Let $g(\lam)$ denote the probability when the number of balls is
Poisson with mean $\lam$.  Corollary \ref{t10} gives that
$g(\lam)\sim 2\eps$ for any $\lam$ for which $\lam = L(1+\eps+o(\eps))$.  But
$L\eps \gg \sqrt{L}$ as $L\gg \eps^{-2}$.  Thus in this
range $f(m)\sim g(m)$ by Lemma \ref{poissonfixed}.
\qed

\begin{prop}\label{t25} Let $\eps\ra 0^+$ and $L\gg \eps^{-2}$.
Let $m=L(1-\eps+o(\eps))$.  Let $p_1^{\sR},\ldots,p_{\sL}^{\sR}$ have all
$p_i^{\sR} = \frac{1}{L} + o(\frac{\eps}{L})$.  Let $f(m)$ denote the
probability of $\ESC_{\sL}^{\sR}$ when precisely $m$ balls are placed in
bins $1,\ldots,L$ according to this distribution.  Then $f(m)\sim \eps$.
\end{prop}

\proof Similar to that of Proposition \ref{t24}, now using Corollary \ref{t11}
instead of Corollary \ref{t10}.
\qed
\vskip0.5cm

\noindent

The following result will be used to show that most placement of balls
which are good on the left and right sides are also good in the middle.
This will be a crucial step in order to show that the probability of
$\TREE$ is asymptotic to the probability of $\ESC_{\sL}\wedge \ESC_{\sL}^{\sR}$.

\begin{prop}\label{middle} Let $M$ balls be placed uniformly in
bins $1,\ldots,M$, let $Z_i$ be the number of balls in bin $i$, and
define a walk by $Y_0=0$, $Y_i=Y_{i-1}+Z_i-1$ for $1\leq i\leq M$.
Set $\MIN$ equal to the minimum of $Y_i$, $0\leq i\leq M$.  Assume
$M,s\ra\infty$.  Then
    \beq\label{qq6} \Pr[\MIN < -s\sqrt{M}] = o(1).
    \eeq
\end{prop}

\proof  The proof makes essential use of Lemma \ref{poissonfixed}.
First suppose all $Z_i\sim \Po(1)$, independent.   As
$s\ra\infty$, $\Pr[Y_{\sM} < -s\sqrt{M}] = o(1)$.  Let $F_i$ be the event
that $Y_i<-s\sqrt{M}$, while $Y_j\geq -s\sqrt{M}$ for all $j<i$. If
$F_i$ occurs, then, by the strong Markov property,
    \eq
    \Pr[Y_{\sM}< -s\sqrt{M}|F_i]
    \geq \Pr[Y_{\sM-i}\leq 0]\geq c,
    \en
where $c>0$ uniformly in $M,i$. Therefore, since the $F_i$ are disjoint and
$\bigvee F_i=\{\MIN < -s\sqrt{M}\}$, we obtain that $\Pr[\bigvee F_i]=o(1)$.
In the terminology of Lemma \ref{poissonfixed},
we have $g(M)=o(1)$ as the total number of balls is Poisson with
mean $M$. Since the event $\{\MIN < -s\sqrt{M}\}$ is decreasing,
we also obtain that $f(M)=o(1)$, where $f(M)=\Pr[\MIN < -s\sqrt{M}]$.
Indeed, in this simple case, this can also be obtained directly
by truncating (\ref{xy7}), and thus noting that
$g(M)\geq f(M)\Pr[\Po(M)\leq M]$,  so that  $f(M)=o(1)$.
\qed

\begin{cor}\label{middlecor} Assume $M,s\ra\infty$ and that
$A,B \geq s\sqrt{M}$.  Let $M+B-A$ balls be placed uniformly
in bins $1,\ldots,M$.  Let $Z_i$ be the number of balls in
bin $i$, and define a walk by $Y_0=A$, $Y_i=Y_{i-1}+Z_i-1$ for
$1\leq i\leq M$ so that $Y_M=B$. Set $\MIN$ equal to the
minimum of $Y_i$, $0\leq i\leq M$.  Then
    \beq\label{qq7} \Pr[\MIN \leq 0] = o(1).
    \eeq
\end{cor}
\proof When $A=B$ this is simply Theorem \ref{middle} with the
walk raised by $A$.
If $A<B$, then ignore the first $B-A$ balls so that now the walk goes from
$A$ to $A$.  If $A>B$, then we add $A-B$ fictitious balls so now the walk
goes from $A$ to $A$ and then we lower the walk by $A-B$ so it goes
from $B$ to $B$.  In both cases we have only increased the probability
that the walk hits zero.  In both cases we have reduced to the $A=B$
case and so (\ref{qq7}) holds.
\qed

\subsection{A simple upper bound on $\Pr[\TREE]$}
\label{sec-hard9}
In this section, we combine Corollaries \ref{t10} and \ref{t11}
to prove the upper bound on $\Pr[\TREE]$ in Theorem \ref{probtreethm}.
To obtain this upper bound, it will be useful to relate the
problem of a fixed number of balls to a Poisson number of balls. This
relation is stated in Proposition \ref{prop-Povsfixed}, and will also
be instrumental in the remainder of the proof of Theorem \ref{probtreethm},
as well as in the proof of Theorem \ref{main2}.

Recall that $M= {k\choose 2} - \sum_j T_j$.
Let $\Pr^*$ be the law where the number $Z_i$ of balls in
bin $i$ is a Poisson random variable with mean $\lambda_i$.
We write
    \eq
    \Lambda=\sum_{i=1}^k \lambda_i.
    \en
The laws of  $\TREE$ under $\P^*$ and $\P$  are related as follows:

\begin{prop}
\label{prop-Povsfixed}For every $\lambda_1, \ldots, \lambda_k$,
and every random variable $X$,
    \eq
    \E[XI[\TREE]]=\frac{\E^*[XI[\TREE]]}{\Pr^*[\Po(\Lambda)=k-1]},
    \en
where the tilt is related to $\lambda_1, \ldots, \lambda_k$ by
    \eq
    p_i= \frac{\lambda_i}{\Lambda}.
    \en
\end{prop}

\proof This result is classical when we note that $\TREE=\TREE\wedge \{\sum_{i=1}^k Z_i=k-1\}$,
and the fact that $\sum_{i=1}^k Z_i=k-1$ has law $\Po(\Lambda)$. Therefore, the claim is
identical to the statement that
    \eq
    E[XI[\TREE]]=E^*\Big[XI[\TREE]\Big|\sum_{i=1}^k Z_i=k-1\Big].
    \en
\qed

\noindent
We  continue by using Proposition \ref{prop-Povsfixed} to prove  a simple bound
for the probability of $\TREE$ which is useful in the course of the proof:

\begin{prop}
\label{prop-probtreeub}
If $L=o(k)$, then
\eq
\P[\TREE]\leq (1+o(1)) \Pr~\!\!^*[\ESC_{\sL}]\Pr~\!\!^*[\ESC_{\sL}^{\sR}].
\en
\end{prop}

\proof We use Proposition \ref{prop-Povsfixed} with $X=1$,
    \eq
    \Pr[\TREE]=\frac{\Pr~\!\!^*[\TREE]}{\Pr^*[\Po(\Lambda)=k-1]}.
    \en
Let $\mu_{\sL},\mu_{\sR}$ be the expected number
of balls in the first $L$ and the last $L$ bins respectively.
From (\ref{leftdist}--\ref{rightdist}), we obtain that
$\mu_{\sL}=L(1+\eps+o(\eps))$
and $\mu_{\sR}=L(1-\eps+o(\eps))$.
Let $m_{\sL},m_{\sR}$ be the actual number
of balls in the first $L$ and the last $L$ bins respectively.
Then we use that
    \eq
    \Pr~\!\!^*[\TREE]\leq \sum_{A,B} \Pr~\!\!^*[\ESC_{\sL}\wedge \{m_{\sL}=A\}]
    \Pr~\!\!^*[\ESC_{\sL}^{\sR}\wedge \{m_{\sR}=B\}]
    \Pr~\!\!^*[\sum_{i=L+1}^{k-L-1} Z_i=k-1-A-B],
    \en
since we omit the requirements on the middle bins imposed by
$\TREE$.  However, uniformly in $A,B$,
    \eqalign
    \Pr~\!\!^*[\sum_{i=L+1}^{k-L-1} Z_i=k-1-A-B]
    &=\Pr~\!\!^*[\Po(\Lambda-\mu_{\sL}-\mu_{\sR})=k-1-A-B]\\
    &\leq \Pr~\!\!^*\Big[\Po(\Lambda-\mu_{\sL}-\mu_{\sR})=\lfloor \Lambda-\mu_{\sL}-\mu_{\sR}\rfloor\Big]\nn\\
    &\sim \frac{1}{\sqrt{2\pi (\Lambda -\mu_{\sL}-\mu_{\sR})}}.\nn
    \enalign

Since $\Lambda=k+o(\sqrt{k})$ and $\mu_{\sL}+\mu_{\sR}=o(k)$, we have
that
    \eq
    \Lambda -\mu_{\sL}-\mu_{\sR}=k+o(k),
    \en
so that
    \eq
    \Pr~\!\!^*\Big[\Po(\Lambda-\mu_{\sL}-\mu_{\sR})=\lfloor \Lambda-\mu_{\sL}-\mu_{\sR}\rfloor\Big]
    \sim \frac{1}{\sqrt{2\pi k}}\sim \Pr~\!\!^*[\Po(\Lambda)=k-1].
    \en
Performing the sums over $A,B$ gives that
    \eq
    \Pr[\TREE]\leq (1+o(1)) \Pr~\!\!^*[\ESC_{\sL}]\Pr~\!\!^*[\ESC_{\sL}^{\sR}].
    \en
Application of Corollaries \ref{t10} and \ref{t11} completes the proof:
\beq\label{prtreeupper} \Pr[\TREE] \leq (1+o(1))2\eps^2 \eeq
\qed

\subsection{The Hard Case: $\Pr[\TREE]$}
\label{sec-hard7}
In this section, we complete the proof of Theorem \ref{probtreethm}.
By Proposition \ref{prop-probtreeub}, it suffices to prove a lower bound.

We place $k-1$ balls into bins $1,\ldots,k$ with left-tilt $p$ as
given by (\ref{dist}).  Recall that $m_{\sL},m_{\sR}$ are the actual number
of balls in the first $L$ and the last $L$ bins respectively.
Note that $m_{\sL},m_{\sR}$ have Binomial distributions with $k-1$
coin flips (the balls) and probability of success $\frac{\mu_{\sL}}{k-1}$
and $\frac{\mu_{\sR}}{k-1}$ respectively.
With foresight, we fix $\omega$ such that
    \beq\label{omegareq} \omega\sqrt{L}\ll \sqrt{k}
    \quad
    \mbox{ and }
    \quad
    \omega\ll \sqrt{L} \quad
    \mbox{ and }
    \quad
    \omega\ra +\infty.
    \eeq
The assumed bounds in (\ref{Lbounds1}) allow us to find such $\omega$.
We say that placement of balls is {\em normal} if
$|m_{\sL}-\mu_{\sL}|< \omega \sqrt{L}$ and $|m_{\sR}-\mu_{\sR}|< \omega\sqrt{L}$.

We shall naturally refer to a partial placement of balls into the
left and right sides, leaving the placement into the middle bins
undetermined, as normal if it meets the above criteria.
We first prove an extension of Theorem \ref{probtreethm},
which will also be useful in proving Theorem \ref{main2}:

\begin{theorem}\label{main1plus}
With probability $\sim 2\eps^2$, the event
$\TREE$ occurs {\em and} the placement is normal.
Consequently, $\Pr[\TREE]\sim 2\eps^2$.
\end{theorem}

Clearly,  Theorem \ref{probtreethm} is a consequence of Theorem
\ref{main1plus}. We first describe a simple example.
When $p=k^{-1.4}$ and $L=k^{0.85}$, we set
$\omega=k^{0.001}$.  Now the probability of a placement not being
normal is $o(\eps^2)$ and so may be ignored.
We now extend the proof to all $p'$s with $pk^{3/2}\rightarrow \infty$:

%

\proof
Let $\NICE$ denote the event $\ESC_{\sL}\wedge \ESC_{\sL}^R \wedge
\{m_{\sL},m_{\sR}\}~\normal$.
From Propositions \ref{t24}--\ref{t25},
    \eq
    \Pr[\ESC_{\sL}|m_{\sL}=A]\sim 2\eps, \qquad \text{and}\qquad
    \Pr[\ESC_{\sL}^{\sR}|m_{\sR}=B]\sim \eps,
    \en
for every normal $A$ and $B$.
Thus
    \eqalign\label{newaa}
    \Pr[\NICE]
    &=\sum_{A,B~\normal}\Pr[\ESC_{\sL}\wedge \ESC_{\sL}^R|m_{\sL}=A,m_{\sR}=B]
\Pr[m_{\sL}=A,m_{\sR}=B]\\
    &\sim 2\eps^2\Pr[m_{\sL}, m_{\sR}~\normal]\sim 2\eps^2,\nn
    \enalign
where we use that $\Pr[m_{\sL}, m_{\sR}~\normal]\sim 1$.

We effectively need to show that there is ``no middle sag," that such paths
do not usually hit zero somewhere in the middle.  When $p=k^{-1.4}$
and $L=k^{0.85}$ simple Chernoff bounds give that $\Pr[Y_i\leq 0]$
is exceeding small for any middle $i$.  Summing over all middle $i$
the probability that some middle $i$ has $Y_i\leq 0$ is $o(\eps^2)$
and so may be ignored.  However, the argument for all $p$'s with
$pk^{-3/2}$ is surprisingly delicate.  We will show $\Pr[\TREE|\NICE]
=1-o(1)$.  We shall do this in two steps.

We shall first extend the
paths from $L$ to a larger $L'$ defined below and then complete the
path.

Let $L'$ satisfy
\eq \label{Lprimebounds}
k^{-1/2}p^{-1} \ll L' \ll k
\en
and let $\omega'$ satisfy
    \beq\label{omegaprimereq} \omega'\sqrt{L'}\ll \sqrt{k}
    \quad
    \mbox{ and }
    \quad
    \omega'\ll \sqrt{L'} \quad
    \mbox{ and }
    \quad
    \omega'\ra +\infty.
    \eeq

Let $m_{\sL '}$, $m_{\sR '}$ denote the actual number of balls in the first
$L'$ and the last $L'$ bins respectively and let $\mu_{\sL '}$, $\mu_{\sR '}$
be the expected number of such balls.  We say that a placement of balls
is $L'$-normal if $|m_{\sL '}-\mu_{\sL '}| < \omega'\sqrt{L'}$ and
$|m_{\sR '}-\mu_{\sR '}| < \omega'\sqrt{L'}$.

Let $\NICE'$ denote the event $\ESC_{\sL '}\wedge \ESC_{\sL '}^R \wedge
\{m_{\sL '},m_{\sR '}\}~L'-\normal$.
The arguments yielding
(\ref{newaa}) hold for these values.  That is, $\NICE$ and
$\NICE'$ each hold with probability $\sim 2\eps^2$.

Corollaries \ref{t10}--\ref{t11}
give that $\ESC_{\sL}\wedge \ESC_{\sL}^R$ has probability $\sim 2\eps^2$.  Thus
the probability that $\ESC_{\sL}$ and $\ESC_{\sL}^R$ but that $m_{\sL},m_{\sR}$ are
not both normal is $o(\eps^2)$.  For $\NICE'$ to hold and $\NICE$ to
fail these would all need occur.  Hence $\Pr[\NICE'\wedge\NICE]\sim 2\eps^2$.
That is,
\beq\label{niceniceprime}
\Pr[NICE'|NICE] = 1-o(1)
\eeq

Now we want to show $\Pr[\TREE|NICE']=1-o(1)$.  It suffices to show this
conditioning on explicit normal values $m_{\sL '}, m_{\sR '}$.
Set
$A= 1+m_{\sL '}-L'$ and $B=L'-m_{\sR '}$.  We now consider the middle
bins as those not amongst the first or last $L'$ bins.  In the middle
we are placing balls with left-tilt $p$ and considering a walk that
begins at $A$ and ends at $B$.  Our normality assumption and
(\ref{Lprimebounds}) imply
    \eq
    A\sim B\sim L'\eps \gg \sqrt{k}
    \en
We claim with probability $1-o(1)$, the
walk will not hit zero.  Removing the tilt moves balls to the right,
which makes it more likely that the walk does hit zero. Therefore,
it suffices to show this when the balls are placed with uniform probability
in each bin. This is precisely Corollary \ref{middlecor}, where
$M=k-2L'$.  Note that our selection (\ref{Lprimebounds}) of $L'$
has assured $M\sim k$
and $A,B\gg \sqrt{k}$ so that the conditions of the Corollary are met.

We conclude that conditioning on $\NICE'$ and any particular normal
$m_{\sL '},m_{\sR '}$ the event $\TREE$ holds with probability
asymptotic to one.  Hence
\beq \label{niceprimetree}
\Pr[\TREE|\NICE'] = 1-o(1)
\eeq
Combining this with (\ref{niceniceprime}) gives
\beq \label{nicetree}
\Pr[\TREE|\NICE] = 1-o(1)
\eeq
Combined with (\ref{newaa}), $\Pr[\TREE\wedge\NICE]\sim 2\eps^2$, the
first part of Theorem \ref{main1plus}.  The upper bound (\ref{prtreeupper})
completes the proof.

\qed

\subsection{The Hard Case: Asymptotic Gaussian}
\label{sec-hard8}
In this section, we prove Theorem \ref{main2}. This proof relies on the
rewrite in Proposition \ref{prop-Povsfixed}.
We therefore only need to investigate the law of $M$
under the measure $\Pr^*$. For this, we note that
we can rewrite
    \eq
    \lbeq{TZrew}
    \sum_{j=1}^k T_j=\sum_{i=1}^k iZ_i,
    \en
where $Z_i$ is the number of balls placed in the $i^{\rm th}$ bin.
Recall that $Z_i$ is $\Po(\lambda_i)$, where
    \eq
    \lbeq{lambdadef}
    \lambda_i=k\frac{p(1-p)^i}{1-(1-p)^k}
    \en

Define
    \eq
    \lambda_{i,t}=\lambda_i e^{t(i-k/2)},
    \en
and write $\P_{t}^*$ the law of this process when $Z_i$ is $\Po(\lambda_{i,t})$
for all $i=1, \ldots, k$. We also write $\E_t^*$ for the expectation w.r.t.\
$\P_t^*$.  Note that $\E^*=\E^*_0.$ The proposition below gives an
explicit equality for the moment generating function of $M-\E^*[M]$
conditionally on $\TREE$:

\begin{prop}
\label{prop-totid}
The equality
    \eq
    \lbeq{totid}
    \E_{0}^*(e^{-t(M-\E^*[M])}|\TREE)=
    e^{\sum_{i=1}^k \lambda_i\big[e^{t(i-k/2)}-1-t(i-k/2)\big]}
    \frac{\P_{t}^*[\TREE]}{\P_{0}^*[\TREE]}
    \en
    holds.
\end{prop}
\proof
When $\TREE$ holds, then
    \eq
    Z_1+\ldots+Z_k=k-1.
    \en
Therefore, when $\TREE$ holds, and using \refeq{TZrew},
    \eq
    M= {{k}\choose{{2}}}-\sum_{j=1}^k T_j
    =-\sum_{i=1}^k (i-k/2)Z_i.
    \en
Similarly, since
    \eq
    \sum_{i=1}^k E^*[Z_i]=\sum_{i=1}^k \lambda_i=k-1,
    \en
we also have that
    \eq
    E^*[M]= {{k}\choose{{2}}}-\sum_{j=1}^k T_j
    =-\sum_{i=1}^k (i-k/2)E^*[Z_i],
    \en
so that we arrive at the equality that when $\TREE$ holds
    \eq
    M-\E^*[M]=-\sum_{i=1}^k (i-k/2)(Z_i-\E^*[Z_i]).
    \en
Therefore, we can write out
    \eqalign
    &\E_{0}^*(e^{-t(M-\E^*[M])}|\TREE)\\
    &\quad=\frac{1}{\P_{0}^*[\TREE]}e^{-t\sum_{i=1}^k (i-k/2)\lambda_i}
    \sum_{\vec z\in \N^k}e^{t\sum_{i=1}^k (i-k/2)z_i}\prod_{i=1}^k e^{-\lambda_i}
    \frac{\lambda_i^{z_i}}{z_i!}I[\TREE]\nn\\
    &\quad=\frac{1}{\P_{0}^*[\TREE]}e^{-t\sum_{i=1}^k (i-k/2)\lambda_i}
    \sum_{\vec z\in \N^k}\prod_{i=1}^n e^{-\lambda_i}
    \frac{[e^{t(i-k/2)}\lambda_i]^{z_i}}{z_i!}I[\TREE]\nn\\
    &\quad=\frac{1}{\P_{0}^*[\TREE]}e^{-t\sum_{i=1}^k (i-k/2)\lambda_i}
    e^{\sum_{i=1}^k (\lambda_{i,t}-\lambda_{i})}
    \sum_{\vec z\in \N^n}\prod_{i=1}^n e^{-\lambda_{i,t}}
    \frac{\lambda_{i,t}^{z_i}}{z_i!}I[\TREE]\nn\\
    &\quad=e^{\sum_{i=1}^k \lambda_i[e^{t(i-k/2)}-1-t(i-k/2)]}
    \frac{\P_{t}^*[\TREE]}{\P_{0}^*[\TREE]}.\nn
    \enalign
\qed
\vskip0.5cm

\noindent
We now formulate a corollary of Proposition \ref{prop-totid}. It is statement,
we write $\Pr_{t}$ for the measure where the tilt
is
    \eq
    p_{i,t}=\frac{\lambda_{i,t}}{\Lambda_t}, \qquad \text{where}
    \qquad \Lambda_t =\sum_{i=1}^k \lambda_{i,t}.
    \en
\begin{cor}
\label{cor-WCtot}
Let $t_k=tk^{-3/2}$. Then, for every $t\in \mathbb{R}$ fixed,
    \eq
    \lbeq{asstot}
    \Pr~\!\!\!_{t_k}[\TREE]=2\vep^2(1+o(1)).
    \en
Consequently, for $\Pr$ and conditionally on $\TREE$,
the random variable $k^{-3/2}(M-\E^*[M])$ converges weakly to a
normal distribution with variance $\frac 1{12}$.
\end{cor}

We conclude from Corollary \ref{cor-WCtot} that we obtain the
central limit theorem `for free' from the scaling of the probability
of $\TREE$,  which holds for all $t\in \mathbb{R}$ fixed. As a consequence,
we obtain that Theorem \ref{main2} holds. Therefore, we are left to prove
Corollary \ref{cor-WCtot}.

\proof
The equality in \refeq{asstot} follows by Theorem \ref{main1plus}, using the extensions
in Corollaries \ref{t10}--\ref{t11}.
Indeed, we first check the assumptions on $\lambda_{i,t}$.
We note that when $i=o(k)$, then
    \eq
    \lambda_i-\lambda_{i,t/k^{3/2}}=\lambda_i + O(k^{-1/2}).
    \en
Since $\lambda_i=1+\vep+o(\vep)$ for $i=o(k)$,
we therefore obtain that as long as $\vep\gg k^{-1/2}$ and $i=o(k)$,
    \eq
    \lambda_{i,t/k^{3/2}}=1+\vep+o(\vep).
    \en
Similarly, when $k-i=o(k)$, and again
$\vep\gg k^{-1/2}$,
    \eq
    \lambda_{i,t/k^{3/2}}=1-\vep+o(\vep).
    \en
Therefore, the asymptotics of $\lambda_{i,t/k^{3/2}}$ are the same as those for
$\lambda_{i}$, and we obtain from Theorem \ref{main1plus} that
    \eq
    \frac{\Pr~\!\!_{t_k}[\TREE]}{\Pr[\TREE]}\sim 1.
    \en

To prove the asymptotic normality of $k^{-3/2}(M-\E^*[M])$
conditionally on $\TREE$, we start by using Proposition \ref{prop-Povsfixed},
which implies
that
    \eq
    \E[e^{-t_k(M-\E^*[M])}I[\TREE]]=\frac{\E^*\big[e^{-t_k(M-\E^*[M])}I[\TREE]\big]}{\Pr^*[\Po(\Lambda)=k-1]}.
    \en
Then we use Proposition \ref{prop-totid} to obtain that
    \eq
    \E\big[e^{-t_k(M-\E^*[M])}\big|\TREE\big]=e^{\sum_{i=1}^k
    \lambda_i[e^{t_k(i-k/2)}-1-t_k(i-k/2)]}\frac{\Pr~\!\!_{t_k}[\TREE]}{\Pr[\TREE]}
    \frac{\Pr^*[\Po(\Lambda_{t_k})=k-1]}{\Pr^*[\Po(\Lambda)=k-1]}.
    \en
Furthermore, using that $\lambda_i-1 =O(\vep)$ uniformly in $i$,
    \eqalign
    \sum_{i=1}^k \lambda_i[e^{t_k(i-k/2)}-1-t_k(i-k/2)]
    &=
    \frac 12\sum_{i=1}^k \lambda_i (i-k/2)^2 t_k^2
    +O(\sum_{i=1}^k |i-k/2|^3 t_k^3)\\
    &=\frac 12\sum_{i=1}^k (i-k/2)^2 t_k^2
    +O(\sum_{i=1}^k |i-k/2|^3 t_k^3)+O(\vep)\nn\\
    &=\sum_{i=0}^{k/2}(i-k/2)^2 t^2 k^{-3}
    +O(k^{-1/2}+\vep)=\frac {t^2}{24} +O(k^{-1/2}+\vep).
    \nn
    \enalign
Moreover, since
    \eqalign
    \Lambda_{t_k}&=\sum_{i=1}^k \lambda_{i,t_k}=\sum_{i=1}^k \lambda_{i}
    +\sum_{i=1}^k \lambda_{i}[e^{t_k(i-k/2)}-1]\\
    &=k-1+\sum_{i=1}^k \lambda_{i}t_k(i-k/2)
    +O\big(\sum_{i=1}^k \lambda_{i}\big[t_k(i-k/2)\big]^2\big)\nn\\
    &=k-1+\sum_{i=1}^k t_k(i-k/2)+\sum_{i=1}^k (\lambda_{i}-1)t_k(i-k/2)
    +O\big(\sum_{i=1}^k \lambda_{i}\big[t_k(i-k/2)\big]^2\big)=k+o(k^{1/2}),\nn
    \enalign
the local central limit theorem remains valid and we obtain
    \eq
    \frac{\Pr^*[\Po(\Lambda_{t_k})=k-1]}{\Pr^*[\Po(\Lambda)=k-1]}\sim 1.
    \en
We conclude that
    \eq
    \lbeq{leftidconcl}
    \E_{0}\Big(e^{-\frac{t}{k^{3/2}}(M-\E^*[M])}\big|\TREE\Big)\sim
    e^{t^2/24}.
    \en
Since $e^{t^2/24}$ is the moment generating function of a Gaussian
random variable with mean 0 and variance 1/12, this completes the
proof of Corollary \ref{cor-WCtot}.
\qed

\section*{Acknowledgement}
This project was initiated during visits of both authors to
Microsoft Research in July 2004.
We thank Benny Sudakov and Michael Krivelevich for useful discussions
at the start of this project. J.S.\ thanks Nitin Arora for help with
the asymptotics of $C(k,l)$.
The work of RvdH was supported in part by the Netherlands Organisation for
Scientific Research (NWO), and was performed in part while visiting
the Mittag-Leffler Institute in November 2004.

\end{document}